# High-order Space-time Flux Reconstruction Methods for Moving Domain Simulation

Meilin Yu*
Corresponding author: mlyu@umbc.edu
\* University of Maryland Baltimore County (UMBC), Baltimore, MD 21250, USA

**Abstract:** A high-order space-time flux reconstruction (FR) method has been developed to solve conservation laws on moving domains. In the space-time framework, the moving domain simulation is similar to that on a stationary domain, except that the shape of the space-time elements varies with time (and space when a deforming grid is used). The geometric conservation law (GCL) can be automatically satisfied to the level of the numerical resolution of the space-time schemes when the space-time discretization of the governing partial differential equations (PDEs) can resolve the geometric nonlinearity of curvilinear space-time elements. In this study, a space-time tensor product operation is used to construct the FR formulation, and the Gauss-Legendre quadrature points are used as solution points both in space and time. A dual time stepping method is used to solve the resulting space-time system. As has been proved by Huynh [*J Sci Comput* **96**, 51 (2023)], in the temporal direction, the FR scheme with the Gauss-Legendre solution points is equivalent to the so-called DG-Gauss implicit Runge-Kutta (IRK) scheme when the quadrature rule based on the solution points (i.e. quadrature points used in DG) is sufficiently accurate to integrate the space-time curvilinear elements. Specifically, we show that when linear space-time elements are adopted in moving domain simulations, the temporal FR scheme based on Gauss-Legendre solution points can always guarantee its equivalency to IRK DG-Gauss. The conditions, under which the moving domain simulation with the method of lines are consistent with those using the space-time formulation, are also discussed. The new space-time FR method can achieve arbitrarily high-order spatial and temporal accuracy without numerical constraints on the physical time step in moving domain simulations. The temporal superconvergence property for moving domain simulations have been demonstrated.

# 1 Introduction

Numerical simulations of fluid flows over moving domains pose big challenges on the development of high-order computational fluid dynamics (CFD) methods with the arbitrary Lagrangian-Eulerian (ALE) formulation [1]. Much research has been conducted to improve the high-order CFD methods for moving domain simulations [2, 3, 4, 5, 6]. A key focus in these efforts is to enforce the (discrete) geometric conservation law (GCL) in the mapping from the time dependent domain to a fixed reference domain in the ALE method. The space-time formulation provides a consistent numerical treatment of both spatial and temporal discretizations [7], and can resolve GCL automatically to the level of the numerical resolution of the schemes used in the simulation. Although some research has been conducted for the space-time discontinuous Galerkin (DG) method [8, 9, 10, 11, 7, 12, 13] and streamline-upwind/Petrov-Galerkin (SUPG) method [14, 15], more work is still needed to fully explore the numerical potential of the space-time method for moving domain simulations.

We note that an implicit space-time method with right Radau points in the temporal direction has been developed by Huynh [16] for conservation laws. More recently, Huynh [17] showed the equivalence between the DG-type discretization methods of ordinary differential equations (ODEs) and several implicit Runge-Kutta (IRK) methods, including Radau IA, Radau IIA [18], and DG-Gauss [19, 20] with the assistance of the flux reconstruction (FR) concept. Note that FR was originally developed by Huynh [21, 22] for compact spatial differentiation of partial differential equations (PDEs), and is a generalization of various types of discontinuous finite/spectral element methods, such as DG [23, 24, 25], spectral difference [26, 27], and spectral volume [28]. In [29], we developed a nodal space-time FR method based on Gauss-Legendre points to solve hyperbolic conservation laws on stationary meshes. Temporal superconvergence with a convergence rate of $(2k + 1)$ for a degree $k$ (i.e. a nominal $(k + 1)^{\text{th}}$ order of accuracy) polynomial construction in the temporal direction was observed. According to [17], in the temporal direction the space-time FR scheme based on Gauss-Legendre points is equivalent to the IRK DG-Gauss, which also has a superconvergence rate of $(2k + 1)$ for a $(k + 1)$ stage construction. In this work, we further extend the nodal space-time FR method for moving domain simulations.

## 2  Problem Statement

Without loss of generality, the GCL in a mapping between a fixed computational domain $(\tau, \xi, \eta, \zeta)$ and a moving physical domain $(t, x, y, z)$ can be expressed as

$$\frac{\partial |J|}{\partial \tau} + \frac{\partial}{\partial \xi}(|J|\xi_t) + \frac{\partial}{\partial \eta}(|J|\eta_t) + \frac{\partial}{\partial \zeta}(|J|\zeta_t) = 0. \tag{1}$$

Herein, $|J|$ is the determinant of the Jacobian matrix of the mapping, and $\xi_t, \eta_t$ and $\zeta_t$ are related to the grid velocity $\boldsymbol{V_g}$ as

$$\xi_t = -\boldsymbol{V_g} \cdot \nabla \xi, \qquad \eta_t = -\boldsymbol{V_g} \cdot \nabla \eta, \qquad \zeta_t = -\boldsymbol{V_g} \cdot \nabla \zeta.$$

If the method of lines is used to perform the moving domain simulation, the term $\partial(\cdot)/\partial\tau$ and the terms $\partial(\cdot)/\partial\xi$, $\partial(\cdot)/\partial\eta$, and $\partial(\cdot)/\partial\zeta$ are usually not discretized in a spatiotemporally consistent manner. This can introduce GCL errors into the numerical simulation. Even if high-order temporal and spatial discretization is used to handle every term, GCL errors can still be generated due to that the evaluation of the grid velocity $\boldsymbol{V_g}$ and the Jacobian $|J|$ is not consistent with each other. To eliminate such type of errors, two approaches are usually adopted. One way is to solve the Jacobian from the GCL formula using the same scheme that solves the unsteady flows [2]; and the other is to directly substitute the temporal derivative with its spatial equivalence [4]. We note that no matter which approach is used, systematic errors due to the inconsistency between the grid velocity and time-dependent metrics can still exist if these terms are not handled with sufficient accuracy. To reduce the GCL-related errors, small time steps can be used in moving domain simulations. However, this may substantially decrease simulation efficiency when solving complex moving grid problems from industrial applications.

In the space-time method, moving domain simulations essentially become static space-time domain simulations, and GCL is satisfied automatically to the level of the numerical scheme's resolution. The grid velocity used in the method of lines is implicitly embedded as the slope of the space-time element along the temporal direction. In **Section 3**, we briefly introduce the nodal space-time FR method for general stationary and moving domain simulations, and then carry out convergence tests of the new method, especially for moving grids with large deformation, in the first two subsection of **Section 4**. In **Section 4.3**, we discuss the connections between the space-time method and the method of lines for moving domain simulations, and compare their performance.

## 3  Numerical Methods

### 3.1  Nodal Space-time Flux Reconstruction Method on Moving Domains

Consider the conservation form of the hyperbolic conservation laws,

$$\frac{\partial Q}{\partial t} + \nabla \cdot \boldsymbol{F}(Q) = 0, \tag{2}$$

defined on $\Omega \times (0, t_{final}]$ with a moving/deforming spatial domain $\Omega$ bounded by $\partial\Omega$ and the time interval $(0, t_{final}]$, where $Q$ is the vector of conservative variables, and $\boldsymbol{F}$ is the spatial flux including both inviscid and viscous terms. Let $\boldsymbol{x} = (x_1, \cdots, x_d)$ be the spatial coordinates, where the subscript '$d$' stands for the dimension of the problem. Now we introduce the space-time domain $\Omega^{st} = \{(\boldsymbol{x}, t) | 0 < t \leq T, \boldsymbol{x} \in \Omega\}$ for the general stationary and moving domain simulation. An illustration of the curvilinear two-dimensional (2D) and three-dimensional (3D) space-time elements for moving domain simulations is presented in **Figure 1**.

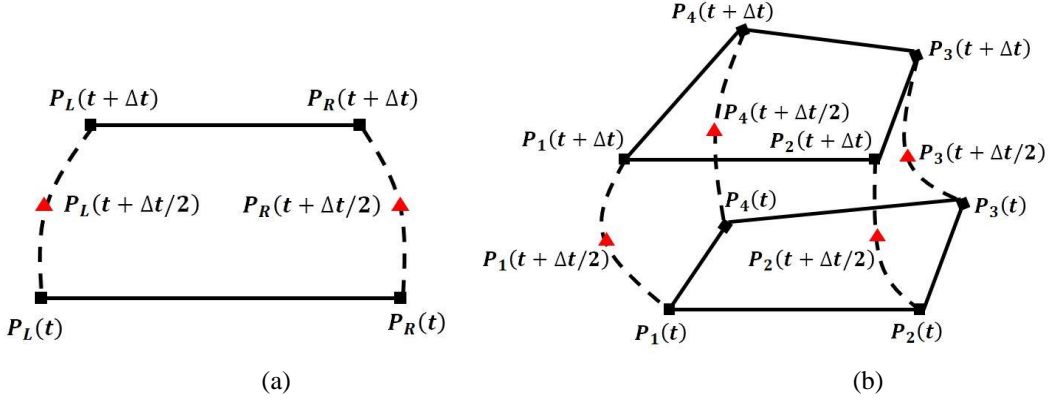

(a)          (b)

**Figure 1**. Curvilinear space-time elements for (a) a 1D moving and deforming spatial (or equivalently 2D spatiotemporal) domain and (b) a 2D moving and deforming spatial (or equivalently 3D spatiotemporal) domain.

On defining the spatiotemporal gradient operator $\nabla_{st} = (\partial_{x_1}, \partial_{x_2}, \cdots, \partial_{x_d}, \partial_t)$, Eq. (2) can then be written as

$$\nabla_{st} \cdot \boldsymbol{F}^{st}(Q) = 0, \tag{3}$$

where $\boldsymbol{F}^{st}(Q) = (\boldsymbol{F}(Q), Q)$ is the space-time flux tensor.

Without loss of generality, we can approximate the exact solution of Eq. (3) using a space-time element-wise continuous polynomial $Q_h(\boldsymbol{x}, t) \in P^k(\Omega^{st})$, where $P^k$ is the polynomial space of order equal to or less than $k$. Let $W(\boldsymbol{x}, t)$ be an arbitrary weighting function or test function. The weighted residual form of the governing equations on each space-time element $\Omega^{st,j}_{t_i \to t_{i+1}}$ (for simplicity and without introducing confusion, this is denoted as $\Omega^{st}$) then reads

$$\int_{\Omega^{st}} \nabla_{st} \cdot \boldsymbol{F}^{st}(Q_h) W dV = 0. \tag{4}$$

Applying integration by parts twice to Eq. (4), one obtains

$$\int_{\Omega^{st}} \nabla_{st} \cdot \boldsymbol{F}^{st}(Q_h) W dV + \int_{\partial\Omega^{st}} \left(F_n^{st,com} - F_n^{st,loc}\right) W dS = 0, \tag{5}$$

where $F_n^{st} = \boldsymbol{n}^{st} \cdot \boldsymbol{F}^{st}$. Note that to ensure conservation, the normal flux term $F_n^{st}$ from the first integration-by-parts operation is replaced with a common flux $F_n^{st,com}(Q_h, Q_h^+, \boldsymbol{n}^{st})$, where $Q_h^+$ denotes the solution outside the current space-time element $\Omega^{st}$, and $\boldsymbol{n}^{st}$ is the outward unit normal vector of $\partial\Omega^{st}$. The normal flux term $F_n^{st}$ from the second integration-by-parts operation is left as its local value $F_n^{st,loc}$ constructed from local solutions. For the inviscid common flux (including $Q$ in the space-time domain) calculation, various approximate Riemann solvers can be used, e.g., Roe approximate Riemann solver [30].

The surface integral in Eq. (5) is then cast into the form of a volume integral via the introduction of a space-time correction field, i.e. $\delta^{st,C} \in P^k(\Omega^{st})$. This is expressed as

$$\int_{\Omega^{st}} \delta^{st,C} W dV = \int_{\partial\Omega^{st}} (F_n^{st,com} - F_n^{st,loc}) W dS. \tag{6}$$

On substituting Eq. (6) into Eq. (5), one obtains

$$\int_{\Omega^{st}} W(\nabla_{st} \cdot \boldsymbol{F}^{st}(Q_h) + \delta^{st,C}) dV = 0. \tag{7}$$

Denote a projection of $\nabla_{st} \cdot \boldsymbol{F}^{st}$ to $P^k(\Omega^{st})$ by $\mathbb{P}(\nabla_{st} \cdot \boldsymbol{F}^{st})$. Then Eq. (7) can be reduced to the differential format as

$$\mathbb{P}(\nabla_{st} \cdot \boldsymbol{F}^{st}) + \delta^{st,C} = 0. \tag{8}$$

This completes a brief derivation of the space-time FR formulation. For an efficient implementation, this formulation can be transformed into a standard (or computational) element. In this study, the Gauss-Legendre points are used in both the spatial and temporal discretization. Note that the temporal discretization is equivalent to the IRK DG-Gauss scheme [17].

### 3.2   Dual Time Stepping Methods

An efficient solution strategy for the nonlinear system resulting from the space-time FR formulation Eq. (8) needs to be developed for good numerical performance. With the dual time stepping procedure, Eq. (8) can be augmented with the pseudo-time derivative term as follows:

$$\frac{\partial Q}{\partial \tilde{t}} + \mathbb{P}(\nabla_{st} \cdot \boldsymbol{F}^{st}(Q_h)) + \delta^{st,C} \rightarrow \frac{\partial Q}{\partial \tilde{t}} = R^{st}(Q_h^n), \tag{9}$$

where $\tilde{t}$ is the pseudo time, and the unsteady residual $R^{st} = -(\mathbb{P}(\nabla_{st} \cdot \boldsymbol{F}^{st}) + \delta^{st,C})$.

In general, the implicit pseudo-time marching can be used to solve Eq. (9), especially when the problem is stiff, such as the flows involving laminar-turbulent transition [31]. Therein, a physics-relevant large time step can be used to speed up numerical simulations. As demonstrated in [32] for the transitional flow over a quasi-3D SD7003 wing at a chord-based Reynolds number of 60,000, the IRK methods, such as the explicit first stage, singly diagonally implicit Runge-Kutta (ESDIRK), can be more cost-effective than explicit time marching. However, for simple unsteady problems, such as those test problems to be used in the next section, an explicit method can cost-effectively drop the relative unsteady residual by 10 orders. Therefore, in this study, the explicit strong stability preserving (SSP) three-stage Runge-Kutta (RK) method [33] is used for pseudo-time marching. It reads

$$\begin{cases} Q^{(1)} = Q^n + \Delta\tilde{t} R^{st}(Q^n) \\ Q^{(2)} = \frac{3}{4}Q^n + \frac{1}{4}Q^{(1)} + \frac{1}{4}\Delta\tilde{t} R^{st}(Q^{(1)}) \\ Q^{n+1} = \frac{1}{3}Q^n + \frac{2}{3}Q^{(2)} + \frac{2}{3}\Delta\tilde{t} R^{st}(Q^{(2)}) \end{cases}. \tag{10}$$

## 4   Results and Discussions

To measure the accuracy, the $L_2$ error of any quantity $u$ within the spatial domain $\Omega$ at a specific time, e.g. the final time $t_{final}$ of the simulation, is used. It is defined as

$$L_2(u_h; \Omega) = \left( \frac{\int_\Omega (u_h^I - u^{exact})^2 dV}{V} \right)^{\frac{1}{2}}. \tag{11}$$

Herein, $u_h$ is the space-time approximation of the quantity $u$, and the upper script '$I$' indicates that the value is an interpolated one calculated from the space-time slab $\Omega \times [t_n - \Delta t, t_n]$ as

$$u_{h,n}^I = \sum_{i=1}^{N} L_i(1) u_{h,i}, \quad L_i(\tau) = \prod_{j=1, j \neq i}^{N} \frac{\tau - \tau_j}{\tau_i - \tau_j}. \tag{12}$$

where $\tau_i \in [-1,1], i = 1, \cdots, N$, are the temporal solution points in the standard space-time element.

### 4.1 Convergence on Stationary Grids

In this subsection, we briefly summarize some numerical experimental results from space-time FR simulations of several 1D and 2D hyperbolic conservation laws on stationary grids [29].

#### 4.1.1 1D Linear Wave Propagation

Consider the 1D linear wave equation,

$$\frac{\partial u}{\partial t} + \frac{\partial f}{\partial x} = 0, \quad f = cu, \text{ and } c > 0, \tag{13}$$

defined on the space-time domain $[0,1] \times (0, +\infty)$ with periodic boundary conditions enforced in space. The initial condition is set as $u_0(x) = sin(2\pi x)$, and $c = 1$ in this test. All simulations were conducted until $t = 1$, during which the wave propagates a distance of one spatial domain size. This propagation time is named one period (1T).

The spatial convergence rates using $P^1$ to $P^5$ spatial constructions, and the temporal convergence rates using $P^1$ to $P^4$ temporal constructions are presented in **Figure 2**. We observe that the spatial convergence rate reaches its optimal value, i.e. the $(k + 1)^{th}$ order is achieved for a degree $k$ polynomial construction; see **Figure 2(a)**. As shown in **Figure 2(b)**, superconvergence with the $(2k + 1)^{th}$ order for a degree $k$ temporal polynomial construction shows up in temporal convergence tests. Since the space-time discretization in the temporal direction with a degree $k$ polynomial construction is equivalent to the IRK DG-Gauss scheme with $(k + 1)$ stages, the order of accuracy of which is $(2k + 1)$, the temporal superconvergence rate of the space-time FR scheme agrees with its theoretical value [17].

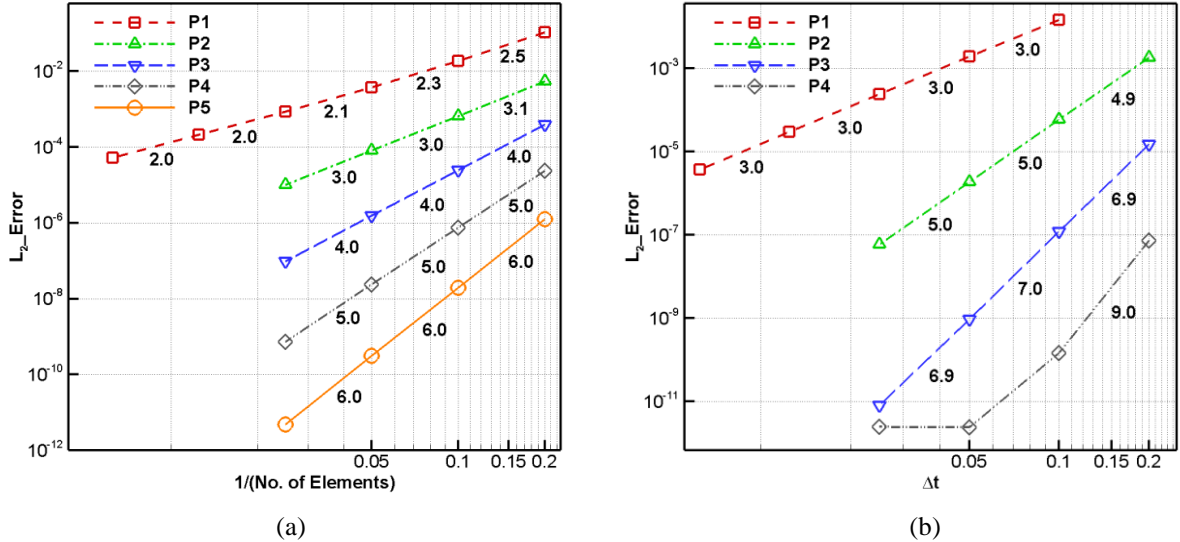

(a) (b)

**Figure 2**. (a) Spatial convergence rates from $P^1$ to $P^5$ spatial constructions for the 1D wave propagation problem on a stationary grid at $t = 1$ (i.e. 1T); and (b) temporal convergence rates from $P^1$ to $P^4$ temporal constructions.

We also measured the $L_2$ error of the quantity $u$ in the space-time slab $\Omega_n^{st} = \Omega \times [t_n - \Delta t, t_n]$ defined as

$$L_2(u_h; \Omega_n^{st}) = \left( \frac{\int_{\Omega_n^{st}} (u_h - u^{exact})^2 dV_{st}}{V_{st}} \right)^{\frac{1}{2}}. \tag{14}$$

From **Figure 3(a)**, we observe that the temporal convergence rate measured from the space-time element is optimal, instead of showing superconvergence, when the simulation was run for 1T. However, for long-term simulations, such as over 500T, superconvergence with the $(2k+1)^{th}$ order for a degree $k$ polynomial construction shows up, as demonstrated in **Figure 3(b)** showing the convergence rate histories of the space-time element based $L_2$ errors for $P^2$ and $P^3$ temporal constructions. Note that for long-term simulations, the spatial convergence rate also presents superconvergent features, i.e. $(2k+1)^{th}$ order for a degree $k$ polynomial construction, as shown in **Figure 4** for $P^2$ and $P^3$ spatial constructions. These observations agree with the superconvergent long-time rate explained in [34].

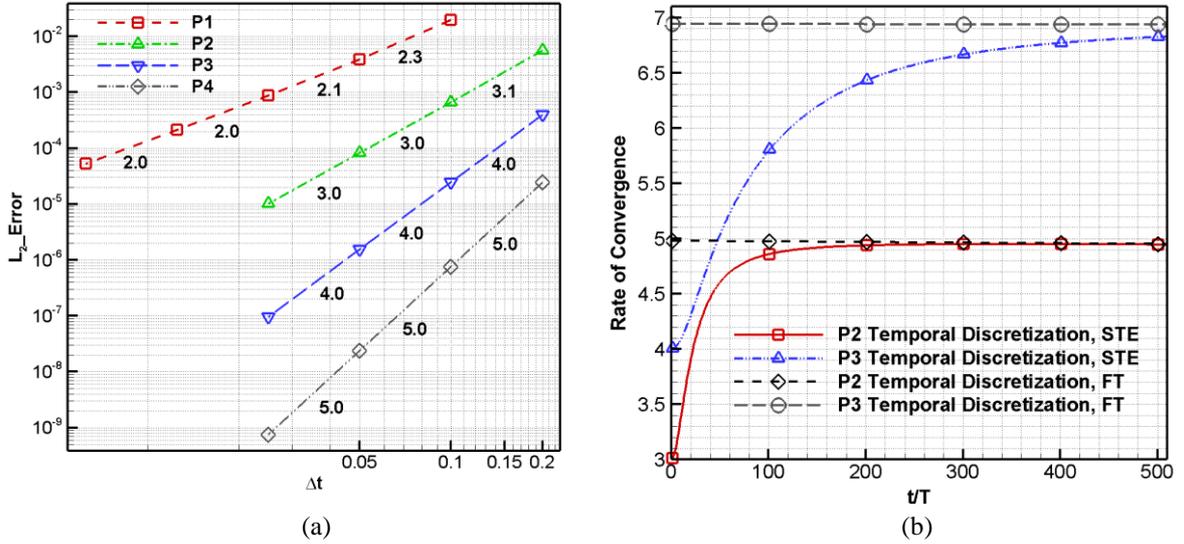

(a) (b)

**Figure 3**. (a) Short-term temporal convergence rates evaluated in the space-time slab $\Omega \times [t_{final} - \Delta t, t_{final}]$ from $P^1$ to $P^4$ temporal constructions for the 1D wave propagation problem on a stationary grid at $t = 1$ (i.e. 1T). (b) Long-term temporal convergence rates evaluated from both the space-time slab $\Omega \times [t_{final} - \Delta t, t_{final}]$ and the final time $t_{final} = 1$ for $P^2$ and $P^3$ temporal constructions.

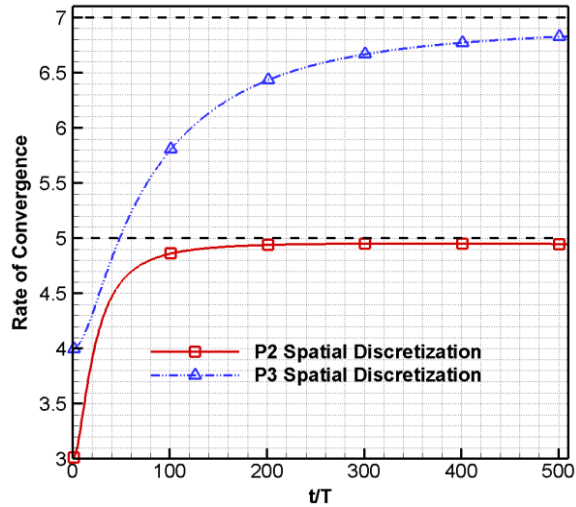

**Figure 4**. Long-term spatial convergence rates for $P^2$ and $P^3$ spatial constructions for the 1D wave propagation problem on a stationary grid.

### 4.1.2 2D Linear Wave Propagation and Euler Vortex Propagation

We further tested the spatial and temporal convergence rates of the space-time FR method for the 2D linear wave propagation and Euler vortex propagation problems. In the 2D wave equation, the spatial flux $\boldsymbol{F}$ has two components $f$ and $g$ with $f = c_1 u$ and $g = c_2 u$, where $c_1$ and $c_2$ are the wave speeds in the $x$- and $y$-directions, respectively. Herein, $c_1 = 0.5$ and $c_2 = 0.5$. In the 2D Euler equation, the conservative variables and fluxes are expressed as

$$Q = \begin{pmatrix} \rho \\ \rho u \\ \rho v \\ \rho E_t \end{pmatrix}, f = \begin{pmatrix} \rho u \\ \rho u^2 + p \\ \rho uv \\ u(\rho E_t + p) \end{pmatrix}, g = \begin{pmatrix} \rho v \\ \rho uv \\ \rho v^2 + p \\ v(\rho E_t + p) \end{pmatrix},$$

where $\rho$ is the fluid density, $u$ and $v$ are velocities in the $x$- and $y$-directions, $p$ is the pressure, and $E_t$ is the specific total energy. The system is closed by the perfect gas law $\rho E_t = p/(\gamma - 1) + \rho(u^2 + v^2)/2$, where $\gamma$ is the specific heat capacity ratio set as 1.4 here. Following [35], one analytical solution for the isentropic vortex can be given as

$$\begin{cases} \rho = \left(1 - \frac{1}{2}(\gamma-1)u_{max}^2 e^{1-\frac{r^2}{b^2}}\right)^{1/(\gamma-1)} \\ p = \frac{1}{\gamma}\left(1 - \frac{1}{2}(\gamma-1)u_{max}^2 e^{1-\frac{r^2}{b^2}}\right)^{\gamma/(\gamma-1)} \\ u = U_0 - \frac{u_{max}}{b} r e^{\frac{1}{2}\left(1-\frac{r^2}{b^2}\right)}\sin\theta \\ v = V_0 + \frac{u_{max}}{b} r e^{\frac{1}{2}\left(1-\frac{r^2}{b^2}\right)}\cos\theta \end{cases}.$$

Herein, $r = \sqrt{(x - U_0 t)^2 + (y - V_0 t)^2}$, $U_0$ and $V_0$ are the advection velocities of the free stream in the $x$- and $y$-directions, respectively, and $\theta$ is the angle with respect to the $x$-direction. In this study, $(U_0, V_0) = (0.5, 0.5)$, $u_{max} = 0.25$, and $b = 0.2$. For both 2D linear wave propagation and Euler vortex propagation problems, the physical domain used in simulations is set as $[-2,2] \times [-2,2]$.

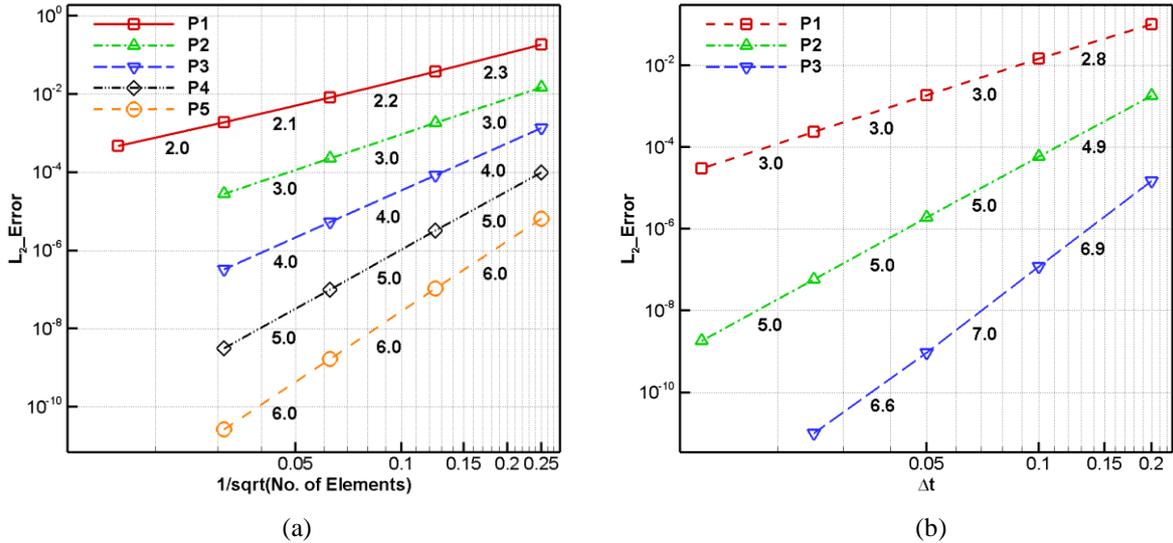

(a)  (b)

**Figure 5**. (a) Spatial convergence rates from $P^1$ to $P^5$ spatial constructions for the 2D wave propagation problem on a stationary grid at $t = 4$ (i.e. 1T); and (b) temporal convergence rates from the $P^1$ to $P^3$ temporal constructions.

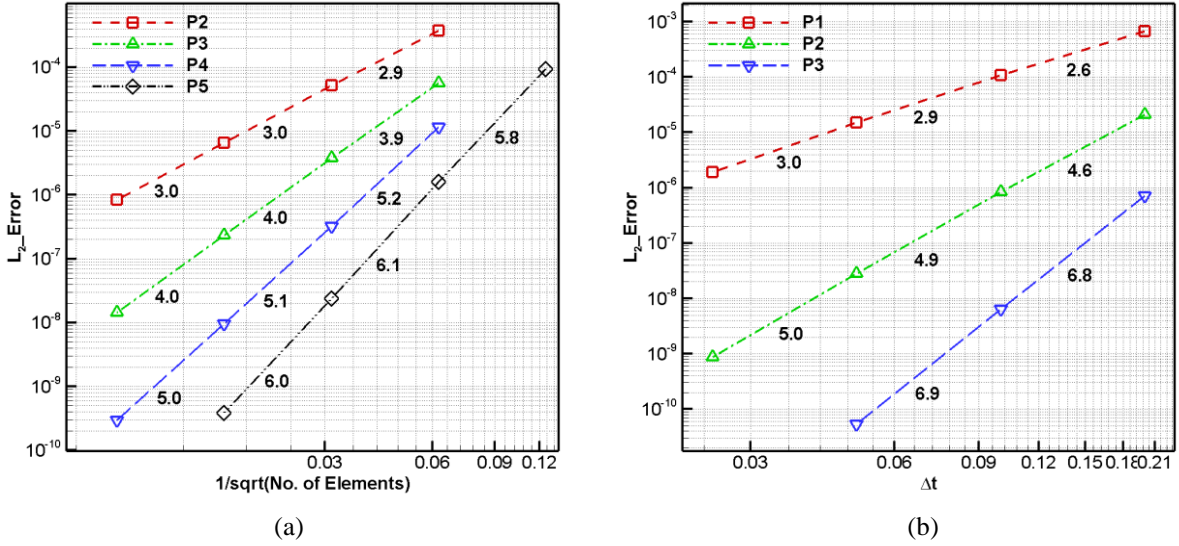

(a)                        (b)

**Figure 6**. (a) Spatial convergence rates from $P^2$ to $P^5$ spatial constructions for the 2D Euler vortex propagation problem on a stationary grid at $t = 0.5$ (i.e. T/8); and (b) temporal convergence rates from the $P^1$ to $P^3$ temporal constructions evaluated at $t = 4$ (i.e. 1T).

For the 2D linear wave propagation, simulations were run until $t = 4$ (i.e. 1T). The spatial convergence rates were tested for $P^1$ to $P^5$ spatial constructions, and the temporal convergence rates were tested for $P^1$ to $P^3$ temporal constructions. We observe from **Figure 5** that the spatial convergence rates match their optimal values, and the temporal convergence rates show the theoretical superconvergent features.

For the 2D Euler vortex propagation, simulations were run until $t = 0.5$ (i.e. T/8) for the spatial convergence tests to save computational costs when simulations were carried out with very high-order schemes on fine grids, and $t = 4$ (i.e. 1T) for the temporal convergence tests. The spatial convergence rates were examined for $P^2$ to $P^5$ spatial constructions, and the temporal convergence rates were tested for $P^1$ to $P^3$ temporal constructions. Similar to the 2D linear wave propagation tests, we observe from **Figure 6** that the spatial convergence rates match their optimal values for short-term simulations, and the temporal convergence rates show the theoretical superconvergent features. This demonstrates that on uniform stationary meshes, spatial nonlinearity from the governing partial differential equations does not affect the temporal superconvergence of the space-time FR schemes.

### 4.2 Convergence on Moving Grids

This section shows the convergence test results for simulations on moving grids. Three types of problems are considered: (1) grid velocity is only a function of time, i.e. there is no grid deformation; (2) deformable grids with spatial and temporal symmetry; and (3) deformable grids without spatial and temporal symmetry.

#### 4.2.1 Moving Grids without Deformation

In this subsection, the nodal space-time FR method is tested on moving grids without deformation. The grid motion is expressed as

$$\begin{cases} x(t) = x_0 + A_x cos(\omega_x t) \\ y(t) = y_0 + A_y cos(\omega_y t) \end{cases}. \quad (15)$$

Herein, $x_0$ and $y_0$ are the initial positions of the coordinates $x$ and $y$ at $t = 0$, $A_x$ and $A_y$ are oscillation amplitudes in the $x$- and $y$-directions, and $\omega_x$ and $\omega_y$ are corresponding oscillation angular frequencies. In this study, the parameters are set as $A_x = A_y = 0.1$, and $\omega_x = \omega_y = 2\pi$. The computational domain is set as $[0,1] \times [0,1]$. All simulations were run until $t = 0.25$ (i.e. T/4) where the grid is located at its maximum offset from its original location.

The spatial convergence rates using $P^2$ to $P^4$ spatial constructions are displayed in **Figure 7(a)**, and the temporal convergence rates using $P^1$ to $P^3$ temporal constructions are shown in **Figure 7(b)**. Since no grid deformation is present, the spatial convergence rates can reach their optimal values, and the temporal convergence rates can reach the theoretical superconvergence of $(2k + 1)^{th}$ order for a $P^k$ temporal polynomial construction. This is similar to the observations from the stationary grids.

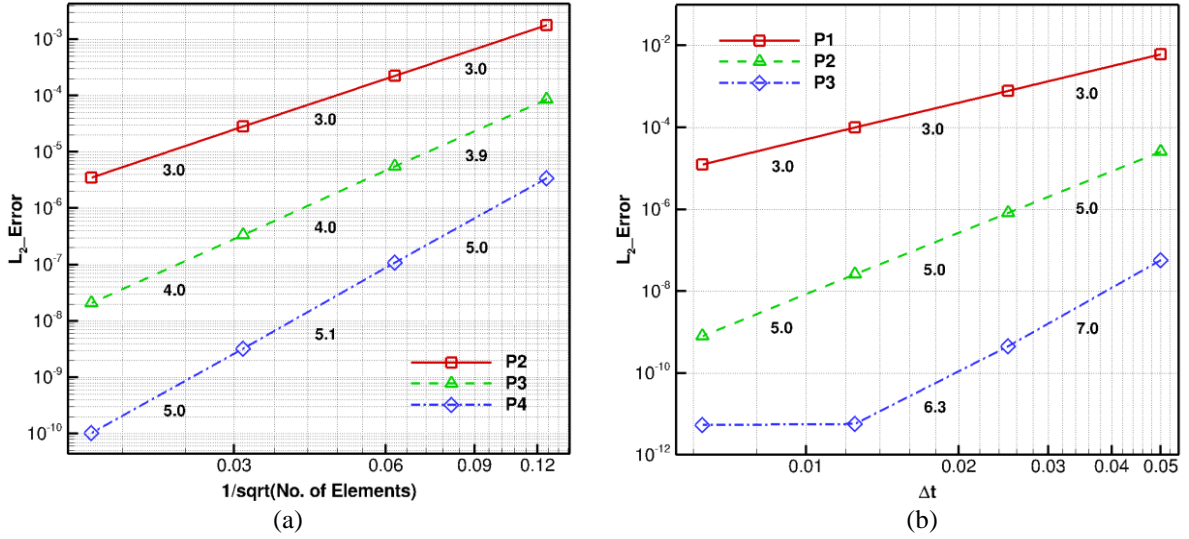

(a)          (b)

**Figure 7**. (a) Spatial convergence rates from $P^2$ to $P^4$ spatial constructions for the 2D linear wave propagation problem on a moving grid without deformation at $t = 0.25$ (i.e. T/4); and (b) temporal convergence rates from $P^1$ to $P^3$ temporal constructions.

### 4.2.2 Moving Grids with Large Grid Deformation

Now we present the convergence test results from simulations on moving grids with large grid deformation. The grid deformation strategy is given as

$$\begin{cases} dx(t) = A_x L_x \, dt/t_{max} \, sin(\omega_t t) sin(\omega_x x) sin(\omega_y y) \\ dy(t) = A_y L_y \, dt/t_{max} \, sin(\omega_t t) sin(\omega_x x) sin(\omega_y y) \end{cases}. \quad (16)$$

Herein, $dx$ and $dy$ are the grid displacement in the *x*- and *y*-directions, $dt$ is the temporal step size, $L_x$, $L_y$ and $t_{max}$ are reference values in the *x*-, *y*- and *t*-directions, and $A_x$ and $A_y$ are scaling factors to control the amplitudes of grid displacement in the *x*- and *y*-directions. The temporal and spatial angular frequencies are defined as

$$\omega_t = n_t \pi/t_{max}, \omega_x = n_x \pi/L_x, \omega_y = n_y \pi/L_y.$$

In this study, the parameters are set as $A_x = A_y = 0.1$, $L_x = L_y = 1.0$, $n_t = 0.5$, $n_x = n_y = 4.0$, and $t_{max} = 0.2$. The computational domain is set as $[0,1] \times [0,1]$, and all simulations were run until $t = t_{max}$.

The wave field on top of the corresponding deformed grid tessellated with $64 \times 64$ elements for the 2D linear wave propagation problem at $t = t_{max}$ is displayed in **Figure 8(a)**. The spatial convergence rates using $P^2$ to $P^4$ spatial constructions are displayed in **Figure 9(a)**. From the grid refinement study, it is observed that the spatial convergence rate for a $P^k$ reconstruction can asymptotically reach its optimal value. We observe that the convergence rate deteriorates on coarse deformable grids. This is due to that linear space-time elements are used in the tests. As shown in **Figure 8(b)**, the $16 \times 16$ grid is plotted on top of the $8 \times 8$ grid. It is clear that the two sets of the grids used in the spatial convergence tests are not conformal to each other. Therefore, the errors caused by the inexact element division in the grid refinement study can contaminate the spatial convergence rates observed over coarse grids. When the grid is refined, this element division related errors becomes smaller, and as a result, the spatial convergence rate asymptotically approaches its optimal values.

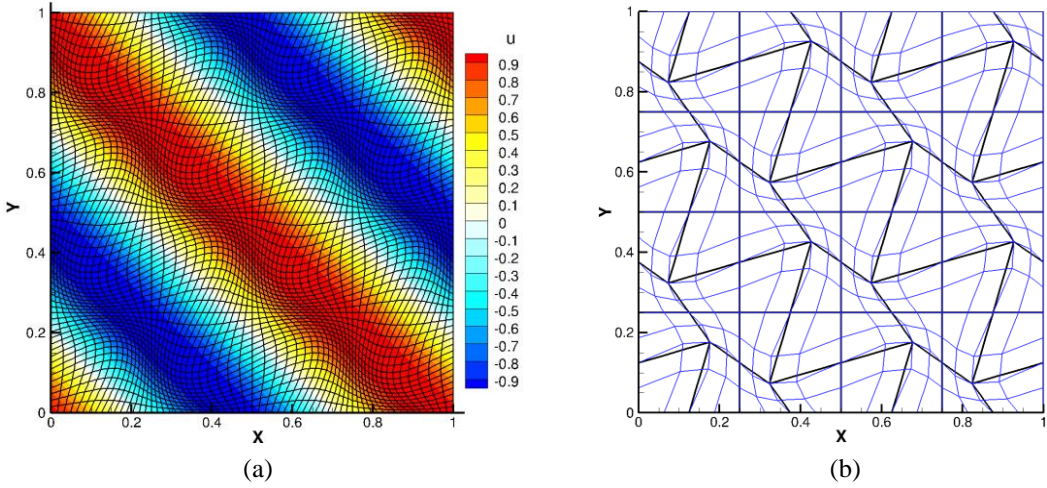

(a)            (b)

**Figure 8**. (a) The flow field and the corresponding grid at $t = t_{max}$ for the 2D linear wave equation on moving grids with large deformation. (b) Two sets of deformed grids with $8 \times 8$ (black) and $16 \times 16$ (blue) elements at $t = t_{max}$.

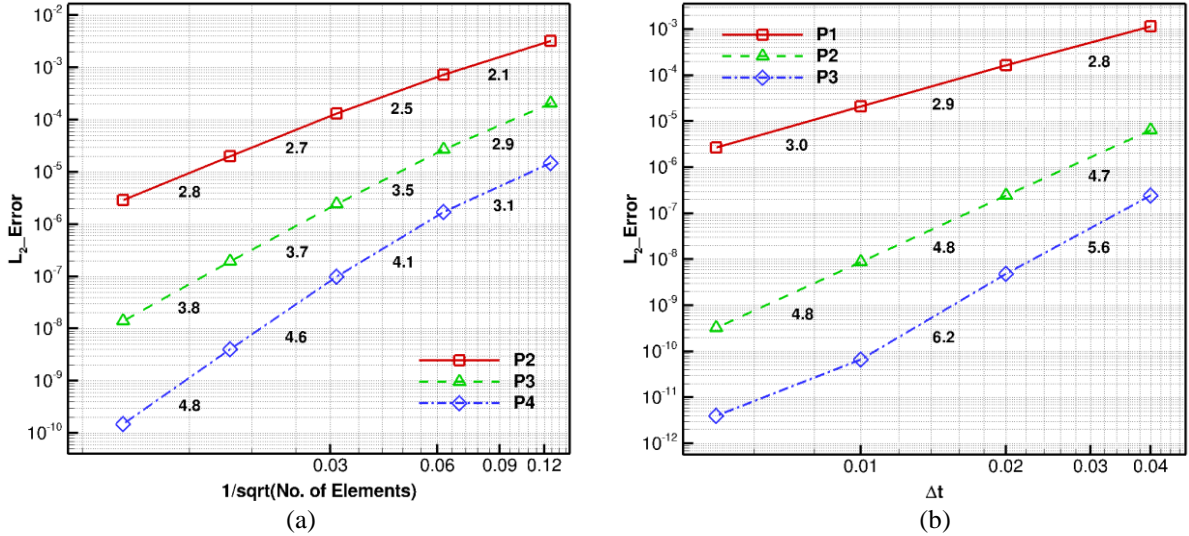

(a)            (b)

**Figure 9**. (a) Spatial convergence rates from $P^2$ to $P^4$ spatial constructions for the 2D linear wave propagation problem on a moving grid with large deformation at $t = t_{max}$; and (b) temporal convergence rates from $P^1$ to $P^3$ temporal constructions.

    The temporal convergence rates using $P^1$ to $P^3$ temporal constructions are displayed in **Figure 9(b)**. Theoretical superconvergence (i.e. $2k + 1$ for a $P^k$ temporal construction) is observed for the $P^1$ and $P^2$ schemes, while the $P^3$ scheme shows a deteriorated order convergence rate. As will be discussed in **Section 4.3**, except for the approxiamtion accuracy of the flow (physics) field, the absolute errors of the space-time methods also depend on the approximation accuracy of the space-time element geometry (i.e. its Jaobian and metrics). Although the linear space-time element can be fully resolved by the space-time FR scheme, it is suspected that the temporal convergence rate deterioration for the $P^3$ scheme is realted to the grid deformation representation capability of linear space-time elements used here. As a side demonstration, the spectral convergence of the temporal discrtization with $\Delta t = 0.01$ is presented in **Figure 10(a)**. Therein, the spectral convergence of the temporal discrtization for the moving grid case without grid deformation (see **Section 4.2.1**), which shows the temporal superconvergence feature, is also presented for comparison. The time step $\Delta t =$

0.05 is used for that case. We note that the temporal error for a nominal $m^{th}$ order (herein, $m = k + 1$) scheme with the superconvergence order $(2m - 1)$ can be approximated as

$$E_m \sim C_m (\Delta t)^{2m-1}, \tag{17}$$

where $C_m$ is a constant depending on the order of accuracy. After taking the logarithm at both sides of Eq. (17), we have

$$\log(E_m) \sim 2 \log(\Delta t)\, m. \tag{18}$$

Therefore, in **Figure 10(a)**, the optimal slope of the curve for the large deformation case with $\Delta t = 0.01$ should be around $2 \log(0.01) = -4$, while that for the no deformation case with $\Delta t = 0.05$ should be around $-2.6$. It is clear that the moving grid case without grid deformation shows the superconvergence feature in the spectral convergence plot; while the moving grid case with large grid deformation fails to maintain superconvergence when the order of accuracy increases. Note that in all tests the order of accuracy of the spatial construction has been maintained sufficiently high (i.e., no apparent error changes are observed beyond the 9$^{th}$ order spatial construction for the $16 \times 16$ grid used in the tests; see **Figure 10(b)**) to ensure that the temporal discretiztion errors dominate. Thus, the resolution constraints on higher-order temporal constructions in the moving grid case with large grid deformation can be relavent to the grid deformation representation capability of linear space-time elements, although the space-time discretization of the governing PDEs can fully resolve the geometric features (i.e. Jabobian and metrics) of linear space-time elements, as will be explained in **Section 4.3**.

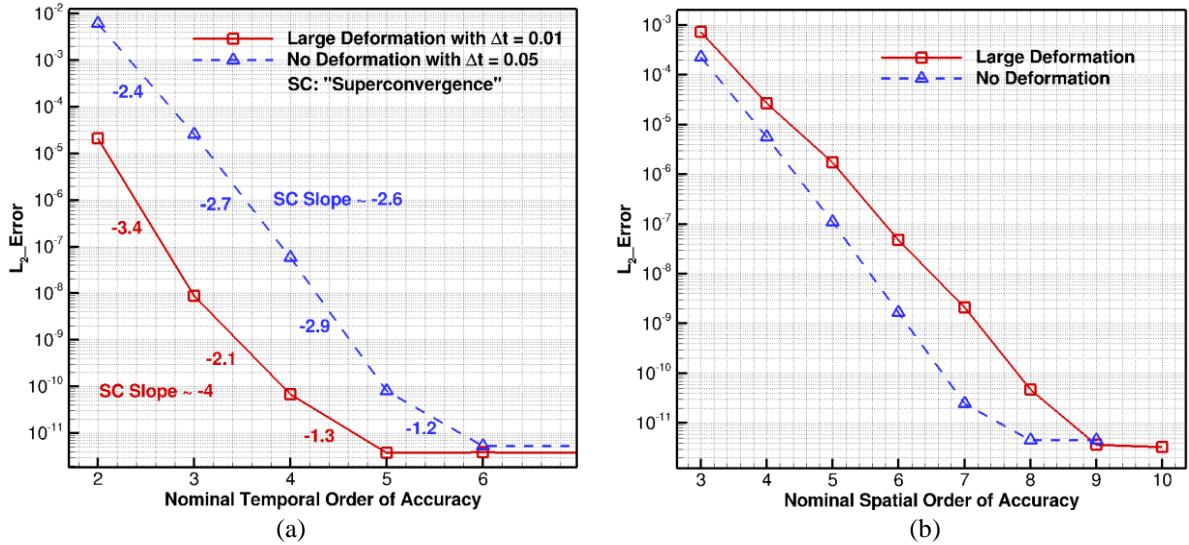

**Figure 10**. (a) Temporal and (b) spatial spectral convergence for the 2D linear wave propagation problem on moving grids at $t = t_{max}$.

### 4.2.3 Tests with a Moving and Deforming Circular Domain

To further study the convergence property of the moving grid space-time FR method, a moving and deforming circular domain case designed in [36] is studied in this subsection. Since the purpose of this study is to test whether the time-varying Jacobian and metrics of a moving grid with spatiotemporal symmetry-breaking large deformation can disturb the flow fields, only the 2D linear wave propagation problem with analytical boundary conditions, instead of the Euler equation with the wall boundary condition in [36], is used for testing purposes. The mesh motion within a circle of a radius 0.5 centered at (0,0) is given as follows:

$$\begin{pmatrix} x(r_0, \theta_0, t) \\ y(r_0, \theta_0, t) \end{pmatrix} = \begin{pmatrix} 1 & 0 & 0 \\ 0 & 1 & \alpha(t) \end{pmatrix} \begin{pmatrix} \cos(A_\theta \alpha(t)) & -\sin(A_\theta \alpha(t)) & 0 \\ \sin(A_\theta \alpha(t)) & \cos(A_\theta \alpha(t)) & 0 \\ 0 & 0 & 1 \end{pmatrix} \tag{19}$$

$$\begin{pmatrix} \psi(t) & 0 & 0 \\ 0 & \dfrac{1}{\psi(t)} & 0 \\ 0 & 0 & 1 \end{pmatrix} \begin{pmatrix} r_0 \cos(\theta_g(r_0,\theta_0,t)) \\ r_0 \sin(\theta_g(r_0,\theta_0,t)) \\ 1 \end{pmatrix}.$$

Herein, $x$ and $y$ are the coordinates of the moving mesh, which depend on the time $t$, and $(r_0, \theta_0)$ is the relative position of the grid point to the origin $(0,0)$ at $t = 0$. The functions $\alpha(t)$, $\psi(t)$, and $\theta_g(r_0, \theta_0, t)$ are defined as follows:

$$\alpha(t) = t^3(8 - 3t)/16,$$
$$\psi(t) = 1 + (A_a - 1)\alpha(t),$$
$$\theta_g(r_0, \theta_0, t) = \theta_0 + A_g f_g(r_0, \theta_0, t),$$

where the function $f_g$ is defined as

$$f_g(r_0, \theta_0, t) = \frac{t^6}{t^6 + 0.01}\left(16 r_0^4 + \eta(t, 10, 0.7)(\cos(32\pi r_0^4) - 1)\right)\eta(\theta_0, 1, 0.7)$$

with the symmetry-breaking perturbation $\eta$ defined as

$$\eta(\lambda, \omega, \tau) = \sin\bigl(\omega\lambda + \tau(1 - \cos(\omega\lambda))\bigr).$$

In the above equations, $A_\theta$ is a constant rotation amplitude, $A_a$ is a constant amplification factor for the deformation of a circle into an ellipse, $A_g$ is a constant volume deformation amplitude. Following [36], they are set as $A_\theta = \pi$, $A_a = 1.5$, and $A_g = 0.15$. In the function $\eta$, the variable $\lambda$ represents the independent variables, such as $t$ and $\theta_0$, the variables $\omega$ and $\tau$ are the angular frequency and phase-lag factor, respectively. The physical domain at $t = 0$ and the deformed one at $t = 1$, together with the corresponding wave fields, are shown in **Figure 11**. We observe that the wave field is not disturbed by the grid motion and deformation.

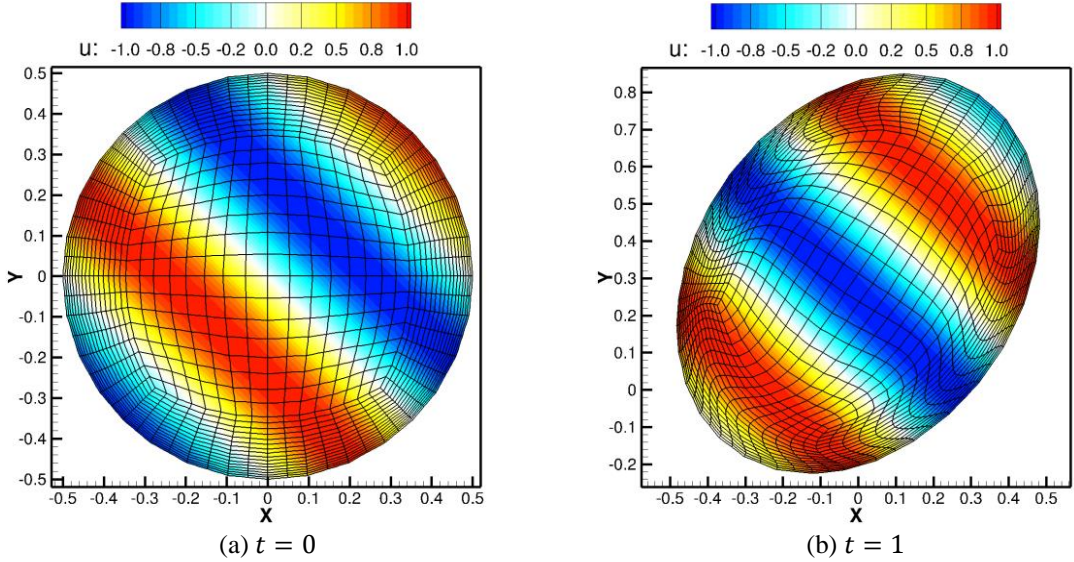

(a) $t = 0$  (b) $t = 1$

**Figure 11**. The flow field and the corresponding grid at (a) $t = 0$ and (b) $t = 1$ for the 2D linear wave equation on a moving and deforming circular domain.

Similar to previous cases, both spatial and temporal convergence rates are examined with a set of grids with 320, 1280, and 5120 elements, respectively. All simulations were run until $t = 1$. $P^2$ to $P^4$ spatial constructions are used to test the spatial convergence rate, and the results are presented in **Figure 12(a)**. It is clear that the optimal short-term convergence rates have been achieved. $P^1$ to $P^3$ temporal constructions are used to test the temporal convergence rate, and the results are presented in **Figure 12(b)**. We observe that the temporal convergence rates for the $P^1$ and $P^2$ schemes are asymptotically approaching their theoretical values, i.e. $2k + 1$, in the tests. However, similar to the case presented in **Section 4.2.2**, the $P^3$ scheme shows the suboptimal superconvergence rate $2k$. Again, it is suspected that the resolution of linear space-time elements causes the convergence rate deterioration. In **Figure 13**, the tragectories $x$ and $y$ of the grid point originally at $(r_0, \theta_0) = (0.4, 1.07274)$ within the time interval $[0,1]$ under different $\Delta t$ conditions are displayed. It shows that although at the end time the $x$ and $y$ are the same for all cases, their tragectories are different. It is not clear yet how these trajectories affect the absolute errors measured with different types of space-time elements, such as linear and curved space-time elements, at the end time, thus the temporal convergence rate. Previous studies [37, 38] have shown that a naïve implementation of curved elements without considering the FR scheme's resolution can cause geometry-related errors. The over-integration treatment used in spectral elements and DG [39, 40] can be a solution to control aliasing errors, and is currently under investigation for moving domain simulations with the space-time FR method.

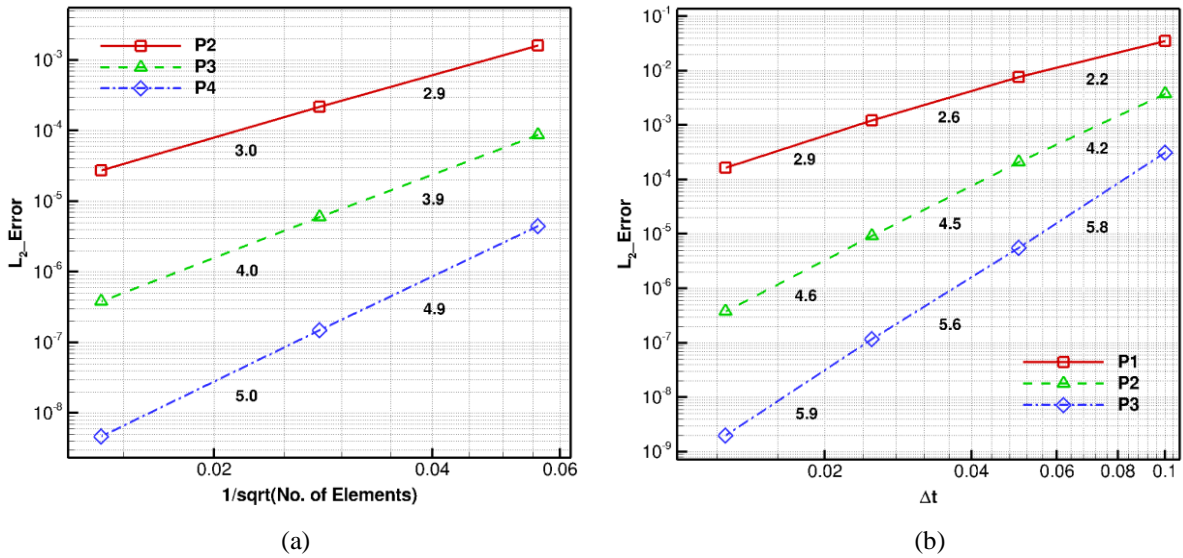

(a)          (b)

**Figure 12**. (a) Spatial convergence rates from $P^2$ to $P^4$ spatial constructions for the 2D linear wave propagation problem on a moving and deforming circular domain at $t = 1$; and (b) temporal convergence rates from $P^1$ to $P^3$ temporal constructions.

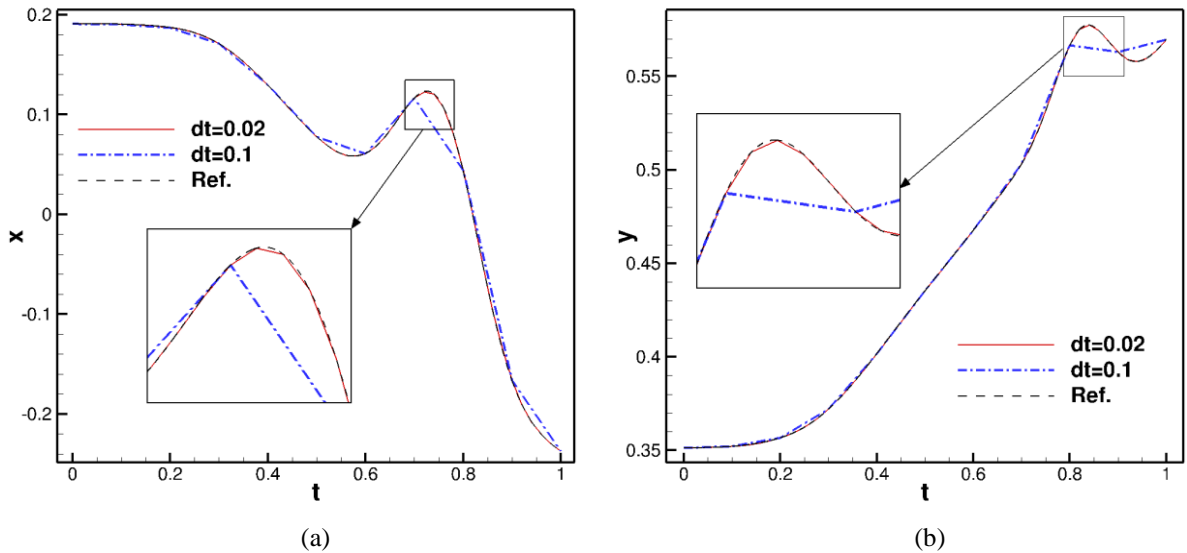

(a)          (b)

**Figure 13**. (a) $x$ and (b) $y$ trajectories of the grid point originally at $(r_0, \theta_0) = (0.4, 1.07274)$ within the time interval [0,1] under different $\Delta t$ conditions.

### 4.3 Comparison with the Method of Lines

We compare the space-time FR method for moving domain simulations to the FR spatial discretization coupled with the method of lines. The purpose is to identify conditions under which the moving domain simulation with the method of lines can be identical to those with the space-time formulation, and to provide the theoretical guidance on how to choose appropriate space-time elements, such as linear and curved space-time elements, for moving domain simulations.

#### 4.3.1 A Space-time Finite Volume View

We use the 1D scalar wave propagation on a 1D moving grid with piece-wise constant velocities at its two ends over the time interval $\Delta t$, and its 2$^{nd}$ order space-time FV discretization as an example to examine some fundamental property of the space-time method for moving domain simulations. As illustrated in **Figure 14**, assume that a 1D spatial element with two ends $x_1$ and $x_2$ at time $t$ has constant velocities at its two ends over the time interval $\Delta t$, and at $t + \Delta t$, the two ends' coordinates become $\tilde{x}_1$ and $\tilde{x}_2$. To facilitate analysis, we always use non-dimensionalized quantities, i.e. the spatial quantity $x$ and time $t$ have been normalized by the characteristic length and time scale, respectively. Therefore, the outward-going normal directions of two side edges of the space-time element shown in **Figure 14** can be expressed as

$$\boldsymbol{n}_1^{st} = \frac{(-\Delta t, \Delta x_1)}{l_1} \text{ and } \boldsymbol{n}_2^{st} = \frac{(\Delta t, -\Delta x_2)}{l_2},$$

where $\Delta x_1 = \tilde{x}_1 - x_1$, $\Delta x_2 = \tilde{x}_2 - x_2$, $l_1 = \sqrt{(\Delta t)^2 + (\Delta x_1)^2}$, and $l_2 = \sqrt{(\Delta t)^2 + (\Delta x_2)^2}$.

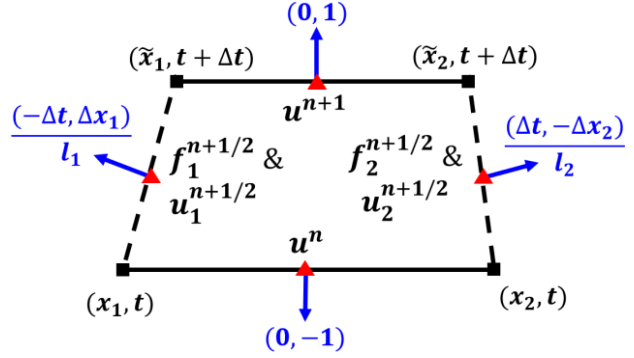

**Figure 14**. An illustration of the 2$^{nd}$ order space-time FV scheme for a 1D moving spatial element.

Next integrate the 1D scalar wave propagation equation in the space-time format over the space-time element shown in **Figure 14**. This reads

$$\int_{\Omega^{st}} \nabla_{st} \cdot \boldsymbol{F}^{st}(u) dV = 0 \rightarrow \int_{\partial \Omega^{st}} \boldsymbol{n}^{st} \cdot \boldsymbol{F}^{st} dS = 0, \tag{20}$$

with $\boldsymbol{F}^{st} = (f(u), u)$. For a 2$^{nd}$ order discretization, the mid-point on each edge is used. Therefore, the discrete form of Eq. (20) is expressed as

$$u^{n+1}(\tilde{x}_2 - \tilde{x}_1) - u^n(x_2 - x_1)$$
$$+ \frac{1}{l_1}\left(-f_1^{n+\frac{1}{2},com}\Delta t + u_1^{n+\frac{1}{2},com}\Delta x_1\right)l_1 + \frac{1}{l_2}\left(f_2^{n+\frac{1}{2},com}\Delta t - u_2^{n+\frac{1}{2},com}\Delta x_2\right)l_2 = 0. \tag{21}$$

Note that the grid velocities at Point 1 and Point 2 are

$$v_{g,1} = \frac{\Delta x_1}{\Delta t} \text{ and } v_{g,2} = \frac{\Delta x_2}{\Delta t}.$$

Therefore, Eq. (21) becomes

$$u^{n+1}(\tilde{x}_2 - \tilde{x}_1) - u^n(x_2 - x_1)$$
$$+ \Delta t\left(-f_1^{n+\frac{1}{2},com} + u_1^{n+\frac{1}{2},com} v_{g,1}\right) + \Delta t\left(f_2^{n+\frac{1}{2},com} - u_2^{n+\frac{1}{2},com} v_{g,2}\right) = 0. \tag{22}$$

Define the 1D normal flux $F = f - u v_g$. Then Eq. (22) reads

$$u^{n+1}(\tilde{x}_2 - \tilde{x}_1) - u^n(x_2 - x_1) + \Delta t \left( F_2^{n+\frac{1}{2},com} - F_1^{n+\frac{1}{2},com} \right) = 0. \tag{23}$$

**Remark 1**. For stationary meshes, $\tilde{x}_2 - \tilde{x}_1 = x_2 - x_1 = \Delta x$, and $v_{g,1} = v_{g,2} = 0$. Therefore, Eq. (23) becomes the Crank–Nicolson scheme if $f^{n+\frac{1}{2}}$ is approximated as $f^{n+\frac{1}{2}} = (f^n + f^{n+1})/2$.

**Remark 2**. Since for a 2$^{nd}$ order FV scheme, the function value at the cell center equals to the cell-averaged value, Eq. (23) can be rewritten as

$$\bar{u}^{n+1}(\tilde{x}_2 - \tilde{x}_1) - \bar{u}^n(x_2 - x_1) + \Delta t \left( F_2^{n+\frac{1}{2},com} - F_1^{n+\frac{1}{2},com} \right) = 0. \tag{24}$$

If a first-order, explicit approximation is used to substitute $F^{n+\frac{1}{2}}$, Eq. (24) finally becomes

$$\bar{u}^{n+1}(\tilde{x}_2 - \tilde{x}_1) - \bar{u}^n(x_2 - x_1) + \Delta t(F_2^{n,com} - F_1^{n,com}) = 0. \tag{25}$$

### 4.3.2  A Finite Volume View from the Method of Lines

Now we derive the 2$^{nd}$ order FV scheme for moving domain simulation in the context of the method of lines. To facilitate understanding, we first transform the governing equation of the 1D scalar wave propagation from the physical domain $[x, t]$ to the computational domain $[\xi, \tau]$, and use $\tau = t$ (therefore $\tau_t = 1$, and $\Delta \tau = \Delta t$) to simplify analysis. Thus, the governing equation in the computational domain can be written as

$$\frac{\partial |J| u}{\partial \tau} + \frac{\partial (|J|\xi_t u + |J|\xi_x f)}{\partial \xi} = 0, \tag{26}$$

with $|J|$ the determinant of the Jacobian matrix of the coordinate transformation, which is $x_\xi$ for the 1D problem, and $|J|\xi_x \equiv 1$ in 1D, and $|J|\xi_t = -x_\tau = -v_g$. The FV formulation can then be written as

$$\int_{\Omega_E} \left( \frac{\partial |J| u}{\partial \tau} + \frac{\partial (|J|\xi_t u + |J|\xi_x f)}{\partial \xi} \right) dV_E = 0, \tag{27}$$

where $\Omega_E = [-1, 1]$ is the fixed standard element in the computational domain. Therefore, we have

$$\int_{\Omega_E} \frac{\partial |J| u}{\partial \tau} dV_E = \frac{\partial}{\partial \tau} \int_{\Omega_E} u |J| \, dV_E = \frac{\partial}{\partial \tau} \int_{\Omega} u \, dV = \frac{\partial}{\partial \tau} (\bar{u} V), \tag{28}$$

and

$$\int_{\Omega_E} \frac{\partial (|J|\xi_t u + |J|\xi_x f)}{\partial \xi} dV_E = \int_{-1}^{1} \frac{\partial (f - u v_g)}{\partial \xi} d\xi. \tag{29}$$

Similar to the space-time analysis, define the 1D normal flux $F = f - u v_g$, and use the nomenclature shown in **Figure 14**. Then Eq. (29) reads

$$\int_{\Omega_E} \frac{\partial (|J|\xi_t u + |J|\xi_x f)}{\partial \xi} dV_E = F_2^{com} - F_1^{com}. \tag{30}$$

On discretizing Eqs. (28) and (30) with the forward Euler scheme, we have

$$\frac{\bar{u}^{n+1} V^{n+1} - \bar{u}^n V^n}{\Delta \tau} + (F_2^{n,com} - F_1^{n,com}) = 0. \tag{31}$$

Based on the 1D moving domain illustrated in **Figure 14**, $V^{n+1} = (\tilde{x}_2 - \tilde{x}_1)$, and $V^n = (x_2 - x_1)$. Moreover, $\Delta \tau = \Delta t$ in the coordination transformation setup. Therefore, Eq. (31) can be written as

$$\bar{u}^{n+1}(\tilde{x}_2 - \tilde{x}_1) - \bar{u}^n(x_2 - x_1) + \Delta t(F_2^{n,com} - F_1^{n,com}) = 0. \tag{32}$$

**Remark 3**. Eq. (32) is the same as Eq. (25) derived from the space-time perspective with the assumption that (1) an at most $2^{nd}$ order FV scheme is used, and (2) forward Euler is used to approximate the fluxes on the spatial element boundaries.

### 4.3.3 Discussion of the FR Formulation with the Method of Lines

Without loss of generality, after extending Eq. (26) to 2D, the governing equation of the scaler wave propagation in the computational domain can be written as

$$\frac{\partial |J|u}{\partial \tau} + \frac{\partial \tilde{f}}{\partial \xi} + \frac{\partial \tilde{g}}{\partial \eta} = 0, \tag{33}$$

where $\tilde{f} = |J|(\xi_t u + \xi_x f + \xi_y g)$ and $\tilde{g} = |J|(\eta_t u + \eta_x f + \eta_y g)$, and $f(u)$ and $g(u)$ are fluxes in the $x$ and $y$ directions, respectively. Following the work [41], we approximate the solution $u$ with the polynomial space sitting on the computational domain $\Omega_E$, i.e., $u_h \in \mathbb{Q}^k(\Omega_E)$. Herein, $\mathbb{Q}^k(\Omega_E) = P^k(\xi) \otimes P^k(\eta)$ is the space of tensor product of polynomials of degree at most $k$ in each variable defined on $\Omega_E$, where $P^k(*)$ is the 1D polynomial space with the polynomial degree at most $k$, and $\otimes$ is the tensor product operator. Note that the determinant $|J|$ of the Jacobian matrix is $|J| = x_\xi y_\eta - x_\eta y_\xi$. Assume that grids tessellated with bilinear elements are used to solve Eq. (33). Thus, $x, y \in \mathbb{Q}^1(\Omega_E)$, and $|J| \in \mathbb{Q}^1(\Omega_E)$. As a result, we have $|J|u_h \in \mathbb{Q}^{k+1}(\Omega_E)$, and Eq. (33) can be expressed in the FR framework as

$$\frac{\partial |J|u_h}{\partial \tau} + \frac{\partial \tilde{f}_h^{loc}}{\partial \xi} + \frac{\partial \tilde{g}_h^{loc}}{\partial \eta} + \frac{\partial \tilde{f}_h^{cor}}{\partial \xi} + \frac{\partial \tilde{g}_h^{cor}}{\partial \eta} = 0, \tag{34}$$

with

$$\begin{cases} \tilde{f}_h = \tilde{f}_h^{loc} + \tilde{f}_h^{cor} \\ \tilde{g}_h = \tilde{g}_h^{loc} + \tilde{g}_h^{cor}, \end{cases} \tag{35}$$

where $\tilde{f}_h^{loc}$ and $\tilde{g}_h^{loc}$ are local fluxes built from the local solution $u_h$, and $\tilde{f}_h^{cor}$ and $\tilde{g}_h^{cor}$ are correction fluxes, which consider the differences between the local fluxes and the common fluxes (also known as numerical fluxes) on element interfaces, including boundary conditions. Define $f_h^{cor} = \tilde{f}_h^{cor}/|J|$ and $g_h^{cor} = \tilde{g}_h^{cor}/|J|$. For quadrilateral elements, the correction fluxes $f_h^{cor}$ and $g_h^{cor}$ are defined as,

$$\begin{aligned} f_h^{cor}(\xi, \eta) &= \Delta f_{h,L}(\eta) G_L(\xi) + \Delta f_{h,R}(\eta) G_R(\xi), \\ g_h^{cor}(\xi, \eta) &= \Delta g_{h,L}(\xi) G_L(\eta) + \Delta g_{h,R}(\xi) G_R(\eta), \end{aligned} \tag{36}$$

where $G_L(*)$ and $G_R(*) \in P^{k+1}(*)$ are the correction functions, which map the differences between numerical fluxes and local fluxes on the boundaries into the entire element, and

$$\Delta f_{h,L}(\eta) = \left(\tilde{f}_{h,L}^{num}(\eta) - \tilde{f}_h^{loc}(-1,\eta)\right)/|J|, \qquad \Delta f_{h,R}(\eta) = \left(\tilde{f}_{h,R}^{num}(\eta) - \tilde{f}_h^{loc}(1,\eta)\right)/|J|,$$

$$\Delta g_{h,L}(\xi) = \left(\tilde{g}_{h,L}^{num}(\xi) - \tilde{g}_h^{loc}(\xi,-1)\right)/|J|, \qquad \Delta g_{h,R}(\xi) = \left(\tilde{g}_{h,R}^{num}(\xi) - \tilde{g}_h^{loc}(\xi,1)\right)/|J|.$$

Since $\Delta f_h$ and $\Delta g_{h,L}$ are defined on element edges, they can be approximated with degree $k$ polynomials along their directions.

Now we discuss the approximation degrees of $\tilde{f}_h^{loc}$, $\tilde{g}_h^{loc}$, $\tilde{f}_h^{cor}$ and $\tilde{g}_h^{cor}$. According to the coordinate transformation,

$$|J|\begin{pmatrix} \xi_x & \xi_y & \xi_t \\ \eta_x & \eta_y & \eta_t \\ 0 & 0 & 1 \end{pmatrix} = \begin{pmatrix} y_\eta & -x_\eta & x_\eta y_\tau - y_\eta x_\tau \\ -y_\xi & x_\xi & y_\xi x_\tau - x_\xi y_\tau \\ 0 & 0 & 1 \end{pmatrix},$$

$|J|\xi_x, |J|\xi_y \in P^1(\xi) \otimes P^0(\eta)$, and $|J|\eta_x, |J|\eta_y \in P^0(\xi) \otimes P^1(\eta)$ for bilinear elements. The grid velocity $\boldsymbol{V}_g = (x_\tau, y_\tau)$ relates to the time-dependent metrics $(\xi_t, \eta_t)$ as follows:

$$\begin{cases} \xi_t = -\boldsymbol{V}_g \cdot \nabla \xi \\ \eta_t = -\boldsymbol{V}_g \cdot \nabla \eta \end{cases} \tag{37}$$

Similar to the assumption used in **Figure 14**, let the node velocity of the 2D moving grid be piecewise constant over the time interval $\Delta t$, i.e., $\boldsymbol{V}_g$ is constant over any $\Delta t$. As a result, $|J|\xi_t \in P^1(\xi) \otimes P^0(\eta)$, and $|J|\eta_t \in P^0(\xi) \otimes P^1(\eta)$. Since $\tilde{f}_h^{loc}$ and $\tilde{g}_h^{loc}$ are built locally from $u_h$, we can require that $f_h^{loc}$ and $g_h^{loc} \in \mathbb{Q}^k(\Omega_E)$. As a result, we have $\tilde{f}_h^{loc} \in P^{k+1}(\xi) \otimes P^k(\eta)$ and $\tilde{g}_h^{loc} \in P^k(\xi) \otimes P^{k+1}(\eta)$, and correspondingly, $\partial \tilde{f}/\partial \xi$ and $\partial \tilde{g}/\partial \eta \in \mathbb{Q}^k(\Omega_E)$. According to Eq. (36), $f_h^{cor} \in P^{k+1}(\xi) \otimes P^k(\eta)$ and $g_h^{cor} \in P^k(\xi) \otimes P^{k+1}(\eta)$. Therefore, $\tilde{f}_h^{cor} \in P^{k+2}(\xi) \otimes P^{k+1}(\eta)$ and $\tilde{g}_h^{cor} \in P^{k+1}(\xi) \otimes P^{k+2}(\eta)$, and correspondingly, $\partial \tilde{f}_h^{cor}/\partial \xi$ and $\partial \tilde{g}_h^{cor}/\partial \eta \in \mathbb{Q}^{k+1}(\Omega_E)$. Thus, we can define the correction field $\delta_h^{cor} = (\partial \tilde{f}_h^{cor}/\partial \xi + \partial \tilde{g}_h^{cor}/\partial \eta)/|J| \in \mathbb{Q}^k(\Omega_E)$, or equivalently, $|J|\delta_h^{cor} \in \mathbb{Q}^{k+1}(\Omega_E)$. As a result, Eq. (34) can be written as

$$\frac{\partial |J| u_h}{\partial \tau} + \frac{\partial \tilde{f}_h^{loc}}{\partial \xi} + \frac{\partial \tilde{g}_h^{loc}}{\partial \eta} + |J|\delta_h^{cor} = 0. \tag{38}$$

To find the weak solution of $u_h \in \mathbb{Q}^k(\Omega_E)$, the weighted residual form of Eq. (38) in DG can be written as

$$\int_{\Omega_E} \left( \frac{\partial |J| u_h}{\partial \tau} + \frac{\partial \tilde{f}_h^{loc}}{\partial \xi} + \frac{\partial \tilde{g}_h^{loc}}{\partial \eta} + |J|\delta_h^{cor} \right) w \, dV_E = 0, \tag{39}$$

where $w \in \mathbb{Q}^k(\Omega_E)$ is an arbitrary weight function sitting on the same space as that of $u_h$.

**Remark 4**. Eq. (39) can be treated as a high-order realization of Eq. (23) for the FV scheme.

**Remark 5**. For the bilinear element with nodal piece-wise constant grid velocities over any time interval $\Delta t$, the Gauss-Legendre quadrature with $(k+1)$ points along each dimension of the element, which is accurate for the integration of a polynomial up to degree $(2k+1)$, can accurately integrate Eq. (39).

**Remark 6**. According to the work [38], when Gauss-Legendre points are used as the solution points for FR and as quadrature points for DG, the FR and DG are equivalent to each other. As has been shown in [42] in the context of stationary mixed triangular and quadrilateral elements, the integration of the last term in Eq. (39) has the following property:

$$\int_{\Omega_E} |J| \delta_h^{cor} w \, dV_E = \int_{\Omega_s} \delta_h^{cor} w \, dV = \int_{\partial \Omega_s} (F_n^{com} - F_n^{loc}) w \, dS, \tag{40}$$

where $F_n = (f_h n_x + g_h n_y) - u_h(v_{g,x} n_x + v_{g,y} n_y)$, $\boldsymbol{V}_g \equiv (v_{g,x}, v_{g,y})$ is the grid velocity, and $\boldsymbol{n} = (n_x, n_y)$ is the unit outward-going normal direction along the surface (edge in 2D) $\partial \Omega_s$ of the element $\Omega_s$.

On transforming Eq. (38) back to the physical domain and substituting $\partial |J|/\partial \tau$ with $-(\partial(|J|\xi_t)/\partial \xi + (|J|\eta_t)/\partial \eta)$ from the GCL, we have

$$\frac{\partial u_h}{\partial \tau} + \nabla \cdot \boldsymbol{F}^{loc} - \boldsymbol{V}_g \cdot \nabla u_h + \delta_h^{cor} = 0, \tag{41}$$

where $\boldsymbol{F}^{loc} = (f_h^{loc}, g_h^{loc})$ are the local fluxes. Finally, the method of lines is used to march Eq. (41) in the time direction.

**Remark 7**. Note that the source term $-\boldsymbol{V}_g \cdot \nabla u_h$ in Eq. (41) is from local manipulations of Eq. (38), and will not affect the correction field $\delta_h^{cor}$.

**Remark 8**. According to the work [4], numerical errors from $-(\partial(|J|\xi_t)/\partial\xi + (|J|\eta_t)/\partial\eta)$ only depend on high-order spatial discretization. In our example, since $|J|\xi_t \in P^1(\xi)\otimes P^0(\eta)$, and $|J|\eta_t \in P^0(\xi)\otimes P^1(\eta)$, this term is a constant. Therefore, substituting $\partial|J|/\partial\tau$ with $-(\partial(|J|\xi_t)/\partial\xi + (|J|\eta_t)/\partial\eta)$ avoids violation of the GCL as $|J|$ is not explicitly evolved temporally; instead, the grid velocity will implicitly determine $|J|$ through grid deformation, and at the same time, it will affect the solution $u_h$ consistent with the implicitly determined $|J|$.

### 4.3.4 Convergence Tests

The moving grid case with large grid deformation presented in **Section 4.2.2** is used here to compare the performance of the space-time method and the method of lines. In the space-time framework, a 9$^\text{th}$ order temporal construction is used in the time direction, and the physical time step is fixed at $\Delta t = 0.02$. The linear space-time elements are used in all tests. In the simulations with the method of lines, the explicit 3$^\text{rd}$ order SSP-RK3 method with $\Delta t = 1 \times 10^{-7}$ is used for time marching. The linear spatial elements are used in all tests. The simulation results for the nominal 3$^\text{rd}$ and 4$^\text{th}$ order spatial discretizations are documented in **Table 1** and **Table 2** for comparison.

We note that the convergence of the space-time method is better than that using the method of lines, although when the grid is refined, both schemes are asymptotically approaching the optimal convergence rate. Based on the analysis in the previous two subsections, the explicit Runge-Kutta methods can only asymptotically match the space-time conditions. Therefore, a very small time step needs to be used to ensure that the errors are not affected by the inexact representation of the grid deformation. The small time step causes huge computational cost. Take the 4$^\text{th}$ order spatial discretization on the $64 \times 64$ grid as an example. The space-time method's cost (i.e. CPU hours in this study) is only about $1/33$ of that of the method of lines. More tests are still needed to compare the performance between the space-time method and IRK methods for moving domain simulations.

**Table 1**: Errors and Order of Accuracy for the 3$^\text{rd}$ order Spatial Discretization

| Meshes | Space-Time | | SSP-RK3 | |
|---|---|---|---|---|
| | Errors | Order | Errors | Order |
| $8 \times 8$ | $3.24 \times 10^{-3}$ | - | $3.74 \times 10^{-3}$ | - |
| $16 \times 16$ | $7.34 \times 10^{-4}$ | 2.1 | $1.09 \times 10^{-3}$ | 1.8 |
| $32 \times 32$ | $1.31 \times 10^{-4}$ | 2.5 | $2.21 \times 10^{-4}$ | 2.3 |
| $64 \times 64$ | $2.00 \times 10^{-5}$ | 2.7 | $3.55 \times 10^{-5}$ | 2.6 |

**Table 2**: Errors and Order of Accuracy for the 4$^\text{th}$ order Spatial Discretization

| Meshes | Space-Time | | SSP-RK3 | |
|---|---|---|---|---|
| | Errors | Order | Errors | Order |
| $8 \times 8$ | $2.04 \times 10^{-4}$ | - | $2.46 \times 10^{-4}$ | - |
| $16 \times 16$ | $2.67 \times 10^{-5}$ | 2.9 | $4.35 \times 10^{-5}$ | 2.5 |
| $32 \times 32$ | $2.44 \times 10^{-6}$ | 3.5 | $4.47 \times 10^{-6}$ | 3.3 |
| $64 \times 64$ | $1.93 \times 10^{-7}$ | 3.7 | $3.98 \times 10^{-7}$ | 3.5 |

## 5  Conclusion and Future Work

This work presents the development of high-order space-time FR methods for moving domain simulations, and discusses numerical approaches that can make the method of lines consistent with the space-time method. Specifically, a space-time tensor product operation is used to construct the FR formulation, and the Gauss-Legendre quadrature points are used as solution points both in space and time. As a result, the temporal construction is equivalent to the IRK DG-Gauss scheme when the quadrature rule based on the solution points (i.e. quadrature points used in DG) is sufficiently accurate to integrate the space-time curvilinear elements. By using the linear space-time elements, we demonstrated the temporal superconvergence properties of the space-time FR methods for moving domain simulations. We found that although lower-order temporal constructions, such as $P^1$ and $P^2$ constructions, can show the $(2k+1)^{th}$ order superconvergence on moving grids with large deformation, the superconvergence of higher-order temporal constructions, such as $P^3$ with a theoretically $7^{th}$ order superconvergence, is suboptimal when grids have large deformation. It is suspected that the grid deformation representation capability of linear space-time elements causes the convergence rate deterioration, although linear space-time elements can be fully resolved by the space-time FR scheme. Spatial and temporal spectral convergence properties of the FR schemes have been used to explain this claim. Precise grid motion capturing with high-order curvilinear space-time elements, and the corresponding higher-resolution space-time FR that can resolve curvilinear space-time elements are currently under development.

We also compared the numerical properties of the space-time method and the method of lines for moving domain simulations, and identified the conditions under which the method of lines is consistent with the space-time formulation. We found that the computational cost of an explicit Runge-Kutta method for moving domain simulations can be prohibitively high to maintain a similar accuracy level achieved by the space-time method. However, considering the equivalence between the temporal construction of the space-time FR method and IRK, the space-time method developed in this study may guide the development of IRK, such as the ESDIRK method, to better simulate moving domain problems. This research is also ongoing, and results will be reported later.